\documentclass{gtart_a}
\pdfoutput=1

\title{Universal circles for quasigeodesic flows}

\author{Danny Calegari}
\givenname{Danny}
\surname{Calegari}
\address{Department of Mathematics\\California Institute of 
Technology\\\newline
Pasadena CA, 91125\\USA}
\email{dannyc@its.caltech.edu}
\urladdr{}

\volumenumber{10}
\issuenumber{}
\publicationyear{2006}
\papernumber{51}
\lognumber{0468}
\startpage{2271}
\endpage{2298}

\doi{}
\MR{}
\Zbl{}


\received{15 June 2004}
\revised{10 September 2006}
\published{29 November 2006}
\publishedonline{29 November 2006}
\corresponding{}
\editor{CPR}
\version{}
\accepted{25 October 2006}
\proposed{David Gabai}
\seconded{Benson Farb, Walter Neumann}

\arxivreference{math.GT/0406040}

\AtBeginDocument{\let\bar\wbar\let\tilde\wwtilde\let\hat\what\let\til\wwtilde}

\makeatletter
\def\cnewtheorem#1[#2]#3{\newtheorem{#1}{#3}[section]
\expandafter\let\csname c@#1\endcsname\c@thm}

\newtheorem{thm}{Theorem}[section]
\cnewtheorem{lem}[thm]{Lemma}
\cnewtheorem{cor}[thm]{Corollary}
\cnewtheorem{conj}[thm]{Conjecture}
\cnewtheorem{qn}[thm]{Question}
\cnewtheorem{prob}[thm]{Problem}

\newtheorem{Athm}{Theorem}

\makeautorefname{Athm}{Theorem}

\newtheorem*{universal_circle_theorem}{Theorem A}
\newtheorem*{lamination_theorem}{Theorem B}
\newtheorem*{norm_theorem}{Theorem C}

\theoremstyle{definition}
\cnewtheorem{defn}[thm]{Definition}
\cnewtheorem{construct}[thm]{Construction}
\cnewtheorem{note}[thm]{Notation}

\theoremstyle{remark}
\cnewtheorem{rmk}[thm]{Remark}
\cnewtheorem{exa}[thm]{Example}

\makeatother  

\def\G{\mathscr G}
\def\F{\mathscr F}
\def\E{\mathscr E}
\def\H{\mathbb H}

\def\hol{\mathrm{hol}}
\def\u{\mathrm{univ}}
\def\homeo{\mathrm{Homeo}}
\def\bhomeo{\mathrm{BHomeo}}
\def\inte{\mathrm{int}}

\def\CP{\mathbb{CP}}

\def\bhomeo{\mathrm{BHomeo}}

\begin{document}

\begin{asciiabstract}
We show that if M is a hyperbolic 3-manifold which admits a
quasigeodesic flow, then pi_1(M) acts faithfully on a universal circle
by homeomorphisms, and preserves a pair of invariant laminations of
this circle. As a corollary, we show that the Thurston norm can be
characterized by quasigeodesic flows, thereby generalizing a theorem
of Mosher, and we give the first example of a closed hyperbolic
3-manifold without a quasigeodesic flow, answering a long-standing
question of Thurston.
\end{asciiabstract}

\begin{htmlabstract}
We show that if M is a hyperbolic 3&ndash;manifold which admits
a quasigeodesic flow, then &pi;<sub>1</sub>(M) acts faithfully on a universal circle
by homeomorphisms, and preserves a pair of invariant laminations of this
circle. As a corollary, we show that the Thurston norm can be characterized
by quasigeodesic flows, thereby generalizing a theorem of Mosher, and we
give the first example of a closed
hyperbolic 3&ndash;manifold without a quasigeodesic flow, answering a
long-standing question of Thurston.
\end{htmlabstract}

\begin{abstract}
We show that if $M$ is a hyperbolic $3$--manifold which admits
a quasigeodesic flow, then $\pi_1(M)$ acts faithfully on a universal circle
by homeomorphisms, and preserves a pair of invariant laminations of this
circle. As a corollary, we show that the Thurston norm can be characterized
by quasigeodesic flows, thereby generalizing a theorem of Mosher, and we 
give the first example of a closed
hyperbolic $3$--manifold without a quasigeodesic flow, answering a
long-standing question of Thurston.
\end{abstract}

\maketitle

\section{Introduction}\label{sec1}
\subsection{Motivation and background}
Hyperbolic $3$--manifolds can be studied from many different
perspectives.  A very fruitful perspective is to think of such
manifolds as {\em dynamical objects}.  For example, a very important
class of hyperbolic $3$--manifolds are those arising as mapping tori of
{\em pseudo-Anosov automorphisms of surfaces} (Thurston
\cite{Thurston_II}). Such mapping tori naturally come with a flow, the
suspension flow of the automorphism. In the seminal paper
\cite{Cannon_Thurston}, Cannon and Thurston showed that this
suspension flow can be chosen to be {\em quasigeodesic} and {\em
pseudo-Anosov}. Informally, a flow on a hyperbolic $3$--manifold is
{\em quasigeodesic} if the lifted flowlines in the universal cover are
quasi-geodesics in $\H^3$, and a flow is {\em pseudo-Anosov} if it
looks locally like a semi-branched cover of an Anosov flow. Of course,
such a flow need not be transverse to a foliation by surfaces. See
Thurston \cite{Thurston_II} and Fenley \cite{Fenley_Anosov} for background and
definitions.

Cannon and Thurston use these two properties to
prove that lifts of surface fibers extend continuously to the sphere at infinity
of $\til{M}$, and their boundaries therefore give natural examples of
Peano-like sphere filling curves.

In \cite{Mosher} and \cite{Mosher_flows}, Mosher (following Gabai), generalizes
\cite{Cannon_Thurston} by constructing
examples of pseudo-Anosov flows almost transverse to any finite 
depth co-orientable taut foliation. In \cite{Fenley_Mosher}, Fenley and Mosher
prove that these flows are actually quasigeodesic. In particular, quasigeodesic
flows exist on any closed hyperbolic $3$--manifold which is not a rational
homology sphere.

All known examples of pseudo-Anosov flows are quasigeodesic,
and it is consistent with all known examples that all quasigeodesic 
flows are homotopic (as oriented line fields)
to pseudo-Anosov flows. On the other hand, the existence of a
pseudo-Anosov flow on a $3$--manifold has many important consequences.
For example, a pseudo-Anosov flow gives rise to a pair of
singular stable and unstable foliations, which can be split open to a pair of
transverse {\em genuine laminations}. Genuine laminations were introduced in
Gabai and Oertel \cite{Gabai_Oertel}, and it is known that they certify many useful properties
of a $3$--manifold Gabai and Kazez \cite{Gabai_Kazez2,Gabai_Kazez1}.

Moreover, in Calegari--Dunfield \cite{Calegari_Dunfield}, it is shown
that if $M$ admits a pseudo-Anosov flow $X$, then $\pi_1(M)$ acts
faithfully on a circle by homeomorphisms, and preserves a pair of
laminations of this circle.  Here a {\em lamination of a circle} is
the boundary data inherited by the circle at infinity from a geodesic
lamination of the hyperbolic plane. Again informally, a lamination of
a circle is just a closed, unlinked subset of the space of pairs of
distinct points in $S^1$. We give more precise definitions in the
sequel.

In \cite{Calegari_Dunfield}, these ($1$--dimensional) dynamical properties of
fundamental groups of $3$--manifolds admitting pseudo-Anosov flows were exploited
to give the first example of a hyperbolic $3$--manifold without a pseudo-Anosov
flow. The example is the {\em Weeks manifold}, the smallest known closed
hyperbolic $3$--manifold, which can be obtained for instance by $(\frac 5 1,\frac 5 2)$ surgery
on the Whitehead link in $S^3$.

In this paper, we show that a quasigeodesic flow on a hyperbolic $3$--manifold
gives rise to a similar dynamical package for $\pi_1(M)$: in particular, we
show that for $M$ a hyperbolic $3$--manifold with a quasigeodesic flow, there is
a faithful action of $\pi_1(M)$ on a circle by homeomorphisms, which preserves
a pair of laminations of this circle. In this way, our theory lets us
promote a codimension $2$ structure (a flow on a $3$--manifold) to a
codimension $1$ structure (laminations on a circle). Dually, our
theory reduces the analysis of holonomy on a two-dimensional
leaf space to the dynamics of $\pi_1(M)$ on a compact one-dimensional manifold.
Experience has shown that the theory of actions of groups on
$1$--dimensional objects is rich and profound, whereas the theory of
group actions on $2$--manifolds remains somewhat elusive; therefore we think that
this dimensional reduction is significant.

An important application of our structure theory is that the unit ball
of the (dual) Thurston norm on cohomology can be detected by
quasigeodesic flows. A basic reference for the Thurston norm is
\cite{Thurston_norm}.  A flow on a $3$--manifold is orthogonal to an
oriented $2$--plane distribution. Such a distribution has a
well-defined integral {\em Euler class}, which is an element of
$H^2(M)$. The Milnor--Wood inequality together with our structure
theory implies that the set of classes obtained from quasigeodesic
flows are contained in the unit ball of the dual Thurston
norm. Conversely, Fenley and Mosher \cite{Fenley_Mosher,Mosher_flows}
show that every extremal class is realized by some quasigeodesic
pseudo-Anosov flow. In particular, the unit ball can be characterized
as the convex hull of such flow classes.

As another important corollary, the calculations of
\cite{Calegari_Dunfield} show that the Weeks manifold does not admit a
quasigeodesic flow. This is the first known example of a hyperbolic $3$ manifold
without a quasigeodesic flow, thereby answering in the negative a long-standing
question of Thurston. This example is not sporadic: combined with
an important recent result of Fenley in \cite{Fenley_noaction}, 
our results imply that there are {\em infinitely}
many closed hyperbolic $3$--manifolds without quasigeodesic flows.

\subsection{Orientation convention}

Throughout this paper, we adhere to the convention that all manifolds 
under consideration are orientable. It follows that the leaf space of
all flows are orientable and co-orientable. We remark that it is straightforward to
extend the results in this paper to non-orientable $3$--manifolds.

\subsection{Statement of results}

In \fullref{sec2}, we introduce the class of {\em product covered flows}. These are
flows on a $3$--manifold $M$ which lift to give a product structure $\R \times \R^2$
on the universal cover $\til{M}$, where the flowlines $l$ are the factors
$\R \times \mathrm{point}$. We show that this topological definition is implied
by a geometric condition, that flowlines on $\til{M}$ are {\em uniformly properly
embedded}; that is, that distances in $\til{M}$ as measured along flowlines,
and as measured in the hyperbolic metric, are comparable at every scale. 

In \fullref{sec3} we show that quasigeodesic flows on hyperbolic $3$ manifolds
are {\em uniformly quasigeodesic}. That is, if $X$ is a flow on $M$ such
that flowlines on $\til{M}$ are all quasigeodesic, then there is a uniform
$k$ such that every flowline is actually $k$--quasigeodesic. In particular, such
a flow is uniform in the sense of \fullref{sec2}, and therefore product covered.

It follows that the natural holonomy representation of $\pi_1(M)$ on
the leaf space of the flow on $\til{M}$ is actually an action by homeomorphisms
on the plane $\R^2$.

In \fullref{sec4} we introduce natural equivalence relations on this $\R^2$, and use the
holonomy action of $\pi_1(M)$ to construct two {\em universal circles}, which
parameterize the action at infinity. Our first main result is the following:

\begin{Athm}\label{ucthm}
Let $M$ be a closed oriented hyperbolic $3$--manifold with a quasigeodesic flow $X$.
Then there are faithful homomorphisms 
$$\rho_\u^\pm\co \pi_1(M) \to \homeo^+((S^1_\u)^\pm)$$
where $(S^1_\u)^\pm$ are topological circles,
called the {\em universal circles of $X$}.
\end{Athm}

As remarked in the introduction, this easily leads to the corollary that the
Weeks manifold admits no quasigeodesic flow.

In \fullref{sec5} we compare the two circles $(S^1_\u)^\pm$ and show that they can
be collated into a single {\em master circle} $S^1_\u$. More precisely, we
prove:

\begin{Athm}\label{lamthm}
Let $M$ be a closed oriented hyperbolic $3$--manifold with a quasigeodesic flow $X$. Then
there is a canonical circle $S^1_\u$, a faithful homomorphism
$$\rho_\u\co \pi_1(M) \to \homeo^+(S^1_\u)$$
and natural monotone maps 
$$\phi^\pm\co S^1_\u \to (S^1_\u)^\pm$$
so that: $$\phi^\pm_*(\rho_\u) = \rho_\u^\pm$$
Moreover, there are a pair of laminations $\Lambda^\pm_\u$ of $S^1_\u$
which are preserved by $\pi_1(M)$, and satisfy: 
$$\phi^\pm(\Lambda^\pm_\u) = \Lambda^\pm$$
\end{Athm}

Finally, in \fullref{sec6} we deduce some corollaries of our structure theory.
We prove:

\begin{Athm}
Let $M$ be a closed oriented hyperbolic $3$--manifold. Then the convex hull of
the set of Euler classes $e(X)$ as $X$ varies over the set of quasigeodesic
flows on $M$ is exactly the unit ball of the dual Thurston norm on $H^2(M)$.
\end{Athm}

Note that one direction of this theorem --- that the convex hull contains the
unit ball --- is due to Gabai, Mosher \cite{Mosher_flows} 
and Fenley--Mosher \cite{Fenley_Mosher}.

This theorem generalizes an earlier theorem of Mosher 
(see \cite{Mosher_efficient_intersection} and \cite{Mosher_norm_II}), 
who showed that the Euler class of a {\em pseudo-Anosov} quasigeodesic 
flow is contained in the unit ball of the dual Thurston norm. Mosher
argues that a pseudo-Anosov quasigeodesic flow can be isotoped to meet
any norm-minimizing surface $S$ in a hyperbolic $3$--manifold ``efficiently'';
the relevant inequality follows easily from this. Our methods give a
new and logically independent proof of Mosher's theorem, but do not give
a new proof of the existence of such an efficient isotopy class of a flow.

\subsection{Acknowledgements}

I would like to thank Nathan Dunfield for a number of valuable comments and corrections.
I would especially like to thank Lee Mosher for bringing my attention to some
important references, and for making many detailed and constructive comments
and corrections. While this research was carried out, I was partially
supported by the Sloan foundation.

\section{Product covered flows}\label{sec2}

A {\em flow} on a manifold $M$ is just a $C^1$ action of $\R$. That is, a $C^1$ 
map $$X\co \R \times M \to M$$ 
such that $X(t,X(s,m)) = X(t+s,m)$ for all $t,s \in \R$ and $m \in M$.
A flow is {\em nonsingular} if this action is locally free; ie, for each
$p \in M$ and each $t \in \R$, the derivative
$$\frac d {dt} X(p,t)$$
is nonvanishing. In this case, the orbits of $X$ give an oriented
$1$--dimensional foliation of $M$, which we denote $X_\F$.
Two different nonsingular flows define the same foliation if and only if
they differ by a reparameterization; conversely, a nonsingular
$1$--dimensional oriented foliation defines a flow by parameterizing
each flowline by arclength.

\begin{defn}
Let $\til{X_\F}$ be the pulled back foliation of a flow $X$ on $M$.
The {\em leaf space} of $\til{X_\F}$ is the quotient space of $\til{M}$
by the equivalence relation $p \sim q$ if $p$ and $q$ are contained in
the same flowline of $\til{X_\F}$.
\end{defn}

The leaf space of $\til{X_\F}$ is rarely Hausdorff, but if it is, it is
a simply connected $2$--manifold. The interesting case for our analysis
will be where this leaf space is $\R^2$ rather than $S^2$

\begin{defn}
A flow $X$ on a $3$--manifold $M$ is called {\em product covered}
if the pullback of the foliation $X_\F$ to the universal cover $\til{M}$
is topologically equivalent to the product foliation of $\R^3$ by
vertical lines.
\end{defn}

If $X$ is product covered,  we denote the quotient map from $\til{M}$ to the leaf
space of $\til{X_\F}$ by
$$\pi_X\co \til{M} \to \R^2$$

\begin{exa}
A linear flow on the torus $T^3$ lifts to a linear flow on $\R^3$, and
is therefore product covered.
\end{exa}

\begin{exa}
Let $M$ be a $3$--manifold, and let $\G$ be an {\em $\R$--covered foliation}.
That is, a co-oriented codimension one foliation such that the space of leaves of
the pullback foliation $\til{\G}$ of the universal cover $\til{M}$ is homeomorphic to
$\R$. A transverse flow $X$ is {\em regulating} if every leaf of $\til{\G}$
intersects every flowline of $\til{X_\F}$. Every regulating flow is
product covered. Moreover, every $\R$--covered foliation admits a regulating
flow. See \cite{Calegari_Rcovered} or \cite{Fenley_Rcovered} for details.
\end{exa}

\begin{defn}
A flow $X$ is {\em uniform} if there is a proper monotone increasing function
$f\co \R^+ \to \R^+$ such that for any flowline $l$ of $\til{X_\F}$, and any two
points $p,q \in l$, there is an inequality
$$d_{\til{M}}(p,q) \ge f(d_l(p,q))$$
with respect to any fixed Riemannian metric on $\til{M}$ pulled back from $M$.
\end{defn}

In other words, there is a {\em uniform} comparison between distance as measured
in a flowline, and distance as measured in $\til{M}$. Notice that with this
definition, if $M$ is a circle bundle over $S^2$, and $X$ is the flow which rotates
the circles at some speed, then $X$ is uniform. On the other hand, this case
is exceptional: we show now that if $M$ is {\em not} a circle bundle over $S^2$, any
uniform flow is product covered.

\begin{lem}\label{uniform_implies_product}
Suppose $X$ is a flow on a closed $3$--manifold $M$. 
If $X$ is uniform and $M$ is not a circle bundle over a sphere, then $X$ is product covered.
\end{lem}
\begin{proof}
Suppose the leaf space of
$\til{X_\F}$ is not Hausdorff. Then there are distinct flowlines $l,m$ of
$\til{X_\F}$ containing points $p \in l, q \in m$, and a sequence of flowlines
$l_i$ of $\til{X_\F}$ and points $p_i,q_i \in l_i$ such that $p_i \to p,q_i \to q$.
Since the sequences $p_i$ and $q_i$ are convergent, 
the distance $d_{\til{M}}(p_i,q_i)$
is eventually bounded by $d_{\til{M}}(p,q) + \epsilon$. On the other hand, since
flowlines converge on compact subsets, the distance between $p_i$ and $q_i$ in
$l_i$ goes to infinity. But this violates uniformity. It follows that
for $X$ uniform, the leaf space of $\til{X_\F}$ is Hausdorff.

If the leaf space
$\til{X_\F}$ is Hausdorff, it is an open, simply connected $2$--manifold;
in particular, it is either a plane, in which case $X$ is product covered, or a sphere,
in which case $\til{M}$ fibers over $S^2$, and $M$ is a circle bundle over $S^2$,
as required.
\end{proof}

\section{Quasigeodesic flows are uniformly quasigeodesic}\label{sec3}

In this section we show that a flow on a hyperbolic $3$--manifold, all of whose
flowlines lift in the universal cover to quasigeodesics, is actually {\em uniformly}
quasigeodesic. In particular, such a flow is uniform, and therefore product covered.
This discussion involves the basic elements of the theory of coarse geometry;
a reference is \cite{Gromov}.

\begin{defn}
Let $X,d_X$ and $Y,d_Y$ be metric spaces. A map
$\phi\co X \to Y$ is a {\em $k$ quasi-isometric embedding} if there is a 
constant $k\ge 1$ such that for all $p,q \in X$,
$$\frac 1 k \, d_Y(\phi(p),\phi(q)) - k \le  d_X(p,q) \le
k \, d_Y(\phi(p),\phi(q)) + k$$
If $X,d_X$ is isometric to $\R$, the image of $X$ under a $k$
quasi-isometric embedding is a {\em $k$ quasigeodesic}.

If $k$ is understood or undetermined, we talk about
quasi-isometric embeddings and quasigeodesics respectively.
\end{defn}

\begin{rmk}
One also says that $\phi\co X \to Y$ is a {\em $(k,\epsilon)$ quasi-isometric embedding} if
there are constants $k\ge 1,\epsilon \ge 0$ such that for all $p,q \in X$,
$$\frac 1 k \, d_Y(\phi(p),\phi(q)) - \epsilon \le  d_X(p,q) \le
k \, d_Y(\phi(p),\phi(q)) + \epsilon$$
The two notions are obviously equivalent, at the cost of possibly increasing the
constants. For the sake of notational simplicity, we prefer to work with the
first definition.
\end{rmk}

The following elementary lemma is well-known. A proof is found in
\cite{Kapovich}.

\begin{lem}\label{straighten_estimate}
Let $\gamma$ be a $k$ quasigeodesic segment, ray or line in $\H^n$. Then there
is a constant $C(k)$ depending only on $k$ and a
geodesic $\gamma_g$ with the same (possibly ideal) endpoints, such that $\gamma$ and
$\gamma_g$ are distance at most $C(k)$ apart in the Hausdorff metric.
\end{lem}

\begin{rmk}
The constant $C(k) = 2^8k^4$ suffices in \fullref{straighten_estimate}.
\end{rmk}

\begin{defn}
A $C^1$ flow $X$ on a $3$--manifold is {\em quasigeodesic} if each
leaf $l$ of $\til{X_\F}$ with its induced length metric
is a quasigeodesic in $\til{M}$.
A flow $X$ is {\em uniformly quasigeodesic} if there is a constant
$k$ such that each leaf $l$ of $\til{X_\F}$ is a
$k$ quasigeodesic.
\end{defn}

\begin{rmk}
Note that since $M$ is compact and $X$ is $C^1$, the parameterizations
of leaves $l$ of $\til{X_\F}$ by arclength or by the flow are uniformly
comparable.
\end{rmk}

\begin{rmk}
The $C^1$ condition is superfluous here. All we really need is that
flowlines are rectifiable, and distance along segments in
flowlines varies continuously in the Hausdorff topology.
\end{rmk}

The following definition is due to Gromov:

\begin{defn}[Gromov]
A path $l$ in $\H^n$ is {\em locally $k$ quasigeodesic on the scale $c$} if for all
pairs of points $p,q \in l$ which are distance 
$t \le c$ apart in $l$, the points $p,q$
are distance $\ge t/k - k$ apart in $\H^n$. If the number
$c$ is understood, we just say $l$ is {\em locally $k$ quasigeodesic}.
\end{defn}

The following lemma is the content of \cite[Remark 7.2.B]{Gromov},
applied to paths in $\H^n$:

\begin{lem}[Gromov]\label{Gromov_lemma}
For every $k\ge 1$ there is a universal scale $c(k)$ such that if
$l$ is a path in $\H^n$ which is locally $k$ quasigeodesic on the scale $c(k)$, 
then $l$ is (globally) $2k$ quasigeodesic. 
\end{lem}

In particular, if $l$ is a path in $\H^n$ which is not $2k$ quasigeodesic, there is
a subsegment of $l$ of length $\le c(k)$ which is not $k$ quasigeodesic. The
function $c(k) = 1000k^2$ works for sufficiently large $k$.

\begin{lem}\label{quasigeodesic_is_uniform}
Let $M$ be a closed hyperbolic $3$--manifold. Then every quasigeodesic
flow $X$ on $M$ is uniformly quasigeodesic.
\end{lem}
\begin{proof}
Suppose to the contrary that we can find a sequence $l_i$ 
of flowlines of $\til{X_\F}$ which are $k_i$ quasigeodesic
for some minimal $k_i$, where $k_i \to \infty$. 
By refining the sequence if necessary,
we can assume $k_i > 2^i$. So $l_i$ is not $2^i$ quasigeodesic, and therefore
by \fullref{Gromov_lemma}, it contains a segment $l_i^i$ of length at most
$c(2^i)$ which is not $2^{i-1}$ quasigeodesic. But then $l_i^i$ contains
a subsegment $l_i^{i-1}$ of length at most $c(2^{i-1})$ which is not 
$2^{i-2}$ quasigeodesic, and so on. Continuing inductively,
we find a nested sequence of segments
$$l_i^n \subset l_i^{n+1} \subset \cdots \subset l_i^i$$
where each $l_i^j$ has length at most $c(2^j)$, and is
not $2^{j-1}$ quasigeodesic. Here $n$ is some universal constant (eg $n=100$)
which should be fixed {\em independently of $i$}.
Let $p_i$ be the midpoint of $l_i^n$,
and fix a sequence of elements $\alpha_i \in \pi_1(M)$
such that $\alpha_i(p_i)$ are contained in a fixed fundamental domain of $\til{M}$.
Choose a convergent subsequence, so that $\alpha_i(p_i) \to p$ where $p$ is
contained in some flowline $l$.
Since the flowlines $\alpha_i(l_i)$ converge on compact subsets, it follows that 
the nested segments
$$\alpha_i(l_i^n) \subset \alpha_i(l_i^{n+1}) 
\subset \cdots \subset \alpha_i(l_i^i) \subset \alpha_i(l_i)$$ 
all converge termwise to 
$$l^n \subset l^{n+1} \subset l^{n+2} \subset \cdots \subset l$$
To see this, observe that
each $\alpha_i(l_i^j)$ has length bounded by $c(2^j)$, and contains a point
$\alpha_i(p_i)$ in a fixed fundamental domain. So the family $\alpha_i(l_i^j)$ is
precompact, and some subsequence converges to some $l^j$. 
Since $l^j$ is contained in $l$, it is unique, and therefore we did not actually
need to pass to a subsequence to get convergence.

By construction, for each $k$, each $l_j^k$ is not $2^{k-1}$ quasigeodesic, so
the same is true of $l^k$. It follows that $l$ is not $2^k$ quasigeodesic for any
$k$. But $l$ is a flowline of $X$, so this violates the hypothesis.

This contradiction proves the lemma.
\end{proof}

\begin{rmk}
Gromov's lemma holds for an arbitrary $\delta$--hyperbolic geodesic metric
space, where now we must
insist that our paths are locally $k$ quasigeodesic on the scale $c(k,\delta)$.
It follows that \fullref{quasigeodesic_is_uniform} 
also applies to $\delta$--hyperbolic
$3$--manifolds.
\end{rmk}

\begin{thm}\label{quasi_is_product}
Let $X$ be a quasigeodesic flow on a closed hyperbolic $3$--manifold. Then $X$
is product covered.
\end{thm}
\begin{proof}
By \fullref{quasigeodesic_is_uniform}, the flowlines of $\til{X_\F}$ are
uniformly $k$ quasigeodesic for some $k$. In particular, the flow $X$ is uniform.
But then the theorem follows immediately from \fullref{uniform_implies_product}
\end{proof}

\begin{note}
If $X$ is a quasigeodesic flow on a closed hyperbolic $3$--manifold $M$, 
we denote the leaf space of $\til{X_\F}$ by $P_X$. Note by \fullref{quasi_is_product}
that $P_X$ is homeomorphic to $\R^2$.
\end{note}

\begin{exa}
The following example is found in \cite{Zeghib}.
Let $M$ be a surface bundle over a circle, and let $\G$ be the foliation
by surfaces. Let $X$ be {\em any} transverse flow. Then flowlines of
$\til{X_\F}$ are uniformly quasigeodesic. To see this, let $\alpha$ be
a nonsingular closed $1$--form whose kernel is tangent to $\G$. Then
$\alpha(X)$ is bounded below by some $\epsilon>0$, whereas $|\alpha|$
(with respect to the hyperbolic metric) is bounded above by some $K$.
It follows that flowlines of $\til{X_\F}$
are uniformly $K/\epsilon$ quasigeodesic.
\end{exa}

\begin{exa}
This generalizes the previous example.
A closed $3$--manifold $M$ is said to {\em slither over $S^1$} if the universal cover of $M$
fibers over $S^1$ 
$$s\co \til{M} \to S^1$$
in such a way that $\pi_1(M)$ acts by bundle automorphisms. Note
that the fibers of $s$ are disconnected.
Suppose $M$ slithers over $S^1$, and suppose $X$ is a flow on $M$
such that for each flowline $l$ of $\til{X_\F}$ the map
$s\co l \to S^1$ is a (universal) covering map.
Then by compactness of $M$, there is a universal constant $C$ such that
every segment $\sigma$ of a flowline $l$ of $\til{X_\F}$ of length at least
$C$ maps surjectively to $S^1$ under $s$. On the other hand, also by
compactness, there is an $\epsilon$ such that every arc in $\til{M}$ of
length $\le \epsilon$ does {\em not} map surjectively to $S^1$ under $s$.
It follows that flowlines of $\til{X_\F}$ are uniformly $C/\epsilon$ quasigeodesic.
\end{exa}

\begin{exa}
In \cite{Mosher}, Mosher constructs an example of a hyperbolic
$3$--manifold $M$ containing a quasifuchsian surface $S$, and a flow
$X$ on $M$ such that away from a single closed orbit $\gamma$, every
flowline intersects $S$. It follows that every flowline of $\til{X_\F}$
either crosses lifts of $S$ with definite frequency, or else spends a
lot of time very close to lifts of $\gamma$. In either case
the flowline is quasigeodesic, and therefore flowlines are uniformly
quasigeodesic.
\end{exa}

\section{Construction of the universal circles}\label{sec4}

In this section we construct two universal circles for a quasigeodesic flow $X$
on a hyperbolic $3$--manifold, and show that there are two natural faithful
homomorphisms from $\pi_1(M)$ to $\homeo^+(S^1)$. Throughout this section
and the following one, we make use of some elementary properties of circular
orders and the order topology. 

\subsection{Circular orders and the order topology}\label{circular_subsection}

We collect here relevant definitions and basic facts about circular orders.

\begin{defn}
Let $S$ be a set. A {\em circular order} on $S$ is a cochain
$$\gamma\co S\times S \times S \to \lbrace -1,0,1 \rbrace$$
which takes the value $0$ if and only if two or more of its co-ordinates are equal,
and which satisfies a {\em cocycle} condition:
$$\gamma(a,b,c) - \gamma(a,b,d) + \gamma(a,c,d) - \gamma(b,c,d) = 0$$
for all quadruples $a,b,c,d \in S$.
\end{defn}

At least for sets $S$ of sufficiently small cardinality, 
one can think of a circular order on $S$ as
an embedding of $S$ in an oriented $S^1$ in such a way that a triple
$(a,b,c)$ takes the value $-1,0,1$ if and only if the three points occur in
clockwise, degenerate or anticlockwise order in $S^1$, with respect to
the orientation on $S^1$. We now make this intuition precise.

A circular order on $S$ determines a topology, called the {\em order topology}.
A basis for this topology consists of sets $U_{a,b}$ with $a,b \in S$
defined by
$$U_{a,b} = \lbrace c \in S \; | \; \gamma(a,b,c) = -1 \rbrace$$
which we call {\em open intervals}, and whose
complements $S\backslash U_{a,b}$ are called {\em closed intervals}.

A circularly ordered set is {\em order complete} if any nested infinite
sequence of closed intervals has a non-empty intersection.

Any circularly ordered set $S$ has a canonical order completion $\overline{S}$
in which $S$ is dense. If $\overline{S}$ is separable, it is compact, and
admits a proper continuous embedding $e\co \overline{S} \to S^1$ 
whose image is closed, and
which is unique up to post-composition with an orientation-preserving homeomorphism of $S^1$.
Any group $G$ of automorphisms of $S$ which preserves the circular
order extends to a group of automorphisms of the order completion $\overline{S}$.
Moreover, since the embedding of $\overline{S}$ in $S^1$ is unique up to homeomorphism,
the action of $G$ on $e(\overline{S})$ extends to an action on $S^1$ by
homeomorphisms.

A reference for this material is \cite{Calegari_euler} or \cite{Thurston_FC2}.

\subsection{Endpoint maps}

We first construct natural {\em endpoint maps} from the leaf space of
$\til{X_\F}$ to the Gromov boundary $S^2_\infty$ of $\til{M} = \H^3$.

\begin{construct}\label{endpoint_maps}
Let $X$ be a quasigeodesic flow on a closed hyperbolic $3$--manifold $M$.
By \fullref{quasigeodesic_is_uniform}, the flowlines of $\til{X_\F}$ are
uniformly $k$ quasigeodesic. By \fullref{quasi_is_product}, the leaf space
$P_X$ of $\til{X_\F}$ is homeomorphic to $\R^2$.
In particular, the action of $\pi_1(M)$ on $\til{M}$
induces a holonomy representation
$$\rho_\hol\co \pi_1(M) \to \homeo^+(P_X)$$
Since each flowline $l$ of $\til{X_\F}$ is a quasigeodesic, 
there are two well-defined
{\em endpoint maps}
$$e^\pm \co  P_X \to S^2_\infty$$
where $S^2_\infty$ denotes the sphere at infinity of $\H^3 = \til{M}$. 
Thinking of a flowline $l$ of $\til{X_\F}$ as a point
in the leaf space $P_X$, we define $e^\pm(l)$ to be the
positive and negative endpoints of the unique oriented geodesic in
$\H^3$ which is a finite Hausdorff distance from $l$.
\end{construct}

\begin{lem}\label{endpoint_maps_continuous}
The endpoint maps $e^\pm$ are continuous.
\end{lem}
\begin{proof}
Suppose $l$ is a complete $k$ quasigeodesic in $\H^3$.
Then by \fullref{straighten_estimate},
the geodesic $l_g$ with the same endpoints is contained in the $2^8k^4$
neighborhood of $l$, and vice versa. 
In particular, if $l_i$ is a sequence of flowlines of 
$\til{X_\F}$ which converges on compact
subsets to $l$, then the straightened geodesics $(l_i)_g$ eventually
contain arbitrarily long segments which are contained in the $2^9k^4$
neigbhorhood of $l_g$. If $(l_i)_g,l_g$ are two hyperbolic geodesics which
are $2^9k^4$ close on a segment of length $t$, then they are $2^9k^4e^{-t}$
close on a subsegment of length $t/2$. It follows that the straightened
geodesics $(l_i)_g$ converge to $l_g$ on compact sets, and therefore the
endpoint maps are continuous.
\end{proof}

\begin{lem}\label{endpoint_maps_dense}
The images $e^\pm(P_X)$ are both dense in $S^2_\infty$.
\end{lem}
\begin{proof}
Since $M$ is a closed hyperbolic $3$--manifold, any non-empty $\pi_1(M)$--invariant subset
of $S^2_\infty$ is dense. The lemma follows.
\end{proof}

\begin{lem}\label{no_interior}
Let $\gamma \subset P_X$ be an embedded circle, and $D \subset P_X$ the
region bounded by $\gamma$. Then there is an equality of images
$$e^\pm(\gamma) = e^\pm(D)$$
\end{lem}
\begin{proof}
For concreteness we concentrate on $e^+$. Suppose to the contrary that there is
$p \in e^+(D)$ which is not in the image of $e^+(\gamma)$.

Let $\sigma\co D \to \til{M}$ be a section of $\pi_X$. That is, we suppose that
$$\pi_X \circ \sigma\co D \to D$$ is the identity.
 
Let $S \subset \til{M}$ be the
union of $\sigma(D)$ and the positive rays contained in the
flowlines of $\til{X_\F}$ which emanate from $\sigma(\gamma)$. Then
the positive rays contained in flowlines of $\til{X_\F}$ which emanate from
$\sigma(D)$ limit to points in $e^+(D)$. Orient $S$ so that the
{\em positive} side of $S$ is the side which contains these positive rays.
Since flowlines of $\til{X_\F}$ are uniformly quasigeodesic, the closure
of $S$ in $\til{M}$ is just $S \cup e^+(\gamma)$.
Then for any point $q \in e^+(D)$ which is not in $e^+(\gamma)$, and for
any sequence $q_i \in \H^3$ limiting to $q \in S^2_\infty$, 
the sequence $q_i$ is eventually
contained on the positive side of $S$.

Now, by \fullref{endpoint_maps_dense}, there is some flowline $l$ of
$\til{X_\F}$ with $e^-(l)$ arbitrarily close to $p$. It follows that the
{\em negative} end of $l$ is contained on the {\em positive} side of $S$.
But if any point $r \in l$ is on the positive side of $S$, then the negative ray
contained in $l$ must intersect $\sigma(D)$, and the negative end of $l$ is
contained on the {\em negative} side of $S$. This contradiction proves the lemma.
\end{proof}

\subsection{Separation properties}

We will study the separation properties in $P_X$ of the various point preimages of
$e^\pm$. The basic algebraic tool we use is Alexander duality.
Since these are arbitrary closed sets, singular cohomology is insufficient for our
purposes.

We recall the definition of the Alexander cohomology theory, following
\cite[Section 6.4]{Spanier}.
For a topological space $X$ and a coefficient module $G$, define $C^q(X;G)$ to
be the module of all functions $\varphi$ from $X^{q+1}$ to $G$ with addition
and scalar multiplication defined pointwise. An element $\varphi \in C^q(X)$
is said to be {\em locally zero} if there is a covering of $X$ by open sets such
that $\varphi$ vanishes on any $(q+1)$--tuple which lies in some element of the
covering. Define $\bar{C}^*(X)$ to be the quotient complex of $C^q(X;G)$ by
the locally zero cochains, and denote the resulting cohomology theory by
$\bar{H}^*(X;G)$.

Alexander cohomology has the following property, which the reader can take as
a working definition:

\begin{thm}{\rm\cite[Corollary 6.9.9]{Spanier}}\label{theory_is_limit}\qua
If $A$ is any closed subset of a manifold $X$, then as $U$ varies over neighborhoods
of $A$ in $X$,
$$\lim_\to \lbrace H^*(U;G) \rbrace \cong \bar{H}^*(A;G)$$
\end{thm}

Moreover, we have the following characterization of connectedness:

\begin{thm}{\rm\cite[Corollary 6.5.7]{Spanier}}\label{connected_characterize}\qua
A nonempty space $X$ is connected if and only if
$$G \cong \bar{H}^0(X;G)$$
\end{thm}

We can now deduce a fundamental property of the point preimages under the
maps $e^\pm$.

\begin{lem}\label{preimages_noncompact}
For every point $p \in S^2_\infty$ in the image of $e^+$, every connected
component of $(e^+)^{-1}(p)$ is unbounded, and similarly for $e^-$.
\end{lem}
\begin{proof}
We suppose not and arrive at a contradiction. 
For ease of notation, define $L=(e^+)^{-1}(p)$. Then $L$ is closed.
Let $K$ be a bounded component of $L$. Then $K$ is compact.

If $L$ is compact, then $L$ can be included into the
interior of a compact disk in $P_X$ whose boundary separates $L$ from infinity,
contrary to \fullref{no_interior}.

Otherwise, $L$ is unbounded. We take $\hat{L}$ to be the closure of $L$ 
in $S^2 = P_X \cup \infty$; ie $\hat{L} = L \cup \infty$. Notice that since $K$
is compact, it is not contained in the connected component which contains
the point $\infty$, and therefore $\hat{L}$ is not connected.

By \fullref{connected_characterize}, we have
$\bar{H}^0(\hat{L};\R) \ne \R$. Moreover, by \fullref{theory_is_limit}, there
is an open set $U \subset S^2$ with $H^0(U;\R) \ne \R$, each component of which
intersects $\hat{L}$. It follows that $U$ has at least two components, and therefore
one component $U_1$ is separated from $U_2: = U\backslash U_1$ with $\infty \in U_2$.
Let $K_1$ be a component of $\hat{L}$ which is contained in $U_1$, and let
$\gamma$ be a loop in $S^2\backslash U$ which separates $U_1$ from $U_2$.
Then $\gamma$ separates $K_1$ from infinity, and bounds a disk $D \subset P_X$
which contains $K_1$ in its interior. On the other hand, 
$e^+(\gamma)$ does not contain $p$. This contradicts \fullref{no_interior} and
proves the lemma.
\end{proof}

\begin{cor}\label{complements_are_disks}
For every point $p \in S^2_\infty$ in the image of $e^+$ and every connected
component $K$ of $(e^+)^{-1}(p)$, the connected
components of $P_X\backslash K$ are unbounded disks.
\end{cor}
\begin{proof}
Since $K$ is connected, closed and unbounded, complementary regions are simply-connected,
by \fullref{theory_is_limit} and standard Alexander duality for polyhedra
in manifolds. If $S$ is a bounded complementary
region, and $p \in S$ is arbitrary, then set $q=e^+(p)$, and
let $L$ be the connected component of $(e^+)^{-1}(q)$ containing $p$.
Since $L$ and $K$ are closed and disjoint, they have disjoint open
neighborhoods. It follows that $L$ is contained in some closed subset of $S$.
In particular, this implies that $L$ is bounded,
contrary to \fullref{preimages_noncompact}.
\end{proof}

\begin{construct}\label{endpoint_relation}
Let $e^\pm\co P_X \to S^2_\infty$ be the endpoint maps from
\fullref{endpoint_maps}. These maps define equivalence relations
$\sim^\pm$ on $P_X$, where the equivalence classes of $\sim^+$ are
the {\em connected components} of the preimages of points in $S^2_\infty$
under $e^+$. The relation $\sim^-$ is defined similarly. 

Since the maps $e^\pm$ are natural, the equivalence classes
of $\sim^\pm$ are permuted by $\pi_1(M)$ and the maps $e^\pm$ factor
through the quotient spaces $T^\pm = P_X/\sim^\pm$. Let
\begin{align*}
\pi^\pm\co& P_X \to T^\pm\\
\tag*{\rlap{\hbox{denote the quotient maps, and}}}
\iota^\pm\co& T^\pm \to S^2_\infty
\end{align*}
the induced maps on the factor spaces, so that $\iota^\pm \circ \pi^\pm = e^\pm$.
\end{construct}

\subsection{Ends and circular orders}

Given a space $X$, the {\em ends} of $X$ (denoted $\E(X)$) may be defined informally
as the connected components of the complement of arbitrarily big compact subsets.
More precisely, if $K_i$ is an exhaustion of $X$ by compact subsets, and
$U_{ij}$ are the connected components of $X\backslash K_i$, each $U_{ij}$ includes
into a unique $U_{i'j'}$ for all $i'<i$, thereby defining a directed system. 
The inverse limit of this directed system is the set $\E(X)$.
If $X$ is a manifold, we may replace the phrase {\em connected component} above
by {\em path connected component}, and think of an end as an equivalence
class of properly embedded ray $r$ in $X$, where $r \sim r'$ if the unbounded
components of $r,r' \cap X\backslash C$ are in the same path component of $X\backslash C$,
where $C$ is any compact set.

Cohomologically, given a collection of (simplicial, singular, Alexander, etc.)
cochains $C^*(X)$ with field coefficients (which we supress in our notation),
there is a natural subcomplex $C^*_c(X)$ consisting of cochains
with compact support, and a natural quotient complex $C^*_e(X)$ defined by
the short exact sequence
$$0 \to C^*_c(X) \to C^*(X) \to C^*_e(X) \to 0$$
Then for reasonable spaces (where reasonable depends on the cohomology theory in
question), $H^0_e(X)$ is a vector space with dimension equal to 
the cardinality of $\E(X)$.

For Alexander cohomology, we have the following interpretation of cohomology with
compact support:

\begin{thm}{\rm\cite[Corollary 6.7.12]{Spanier}}\label{compact_theorem}\qua
If $X$ is a locally compact Hausdorff space and $\hat{X}$ is the one-point compactification
of $X$, there is an isomorphism
$$\bar{H}_c^q(X;G) \cong \til{\bar{H}}^q(\hat{X};G)$$
where the tilde denotes reduced Alexander cohomology.
\end{thm}

For $K \subset P_X$ closed and noncompact, the one point compactification
$\hat{K}$ is just $K \cup \infty \subset S^2$.
Moreover, Alexander duality for polyhedra generalizes as follows:

\begin{thm}{\rm\cite[Theorem 6.9.10]{Spanier}}\label{end_duality_theorem}\qua
Let $X$ be an orientable $n$--manifold. For any closed pair $(A,B)$ in $X$ there is
an isomorphism
$$H_q(X-B,X-A;G) \cong \bar{H}_c^{n-q}(A,B;G)$$
\end{thm}

We deduce that
$$\til{H}_{2-q-1}(P_X\backslash K;G) \cong \bar{H}_c^q(K;G)$$
for any closed $K \subset P_X$.
Note for such $K$ that $\bar{H}^0(K) \cong \R$ if and only
if $K$ is connected, and if $K$ is connected, 
$\bar{H}^0_c(K) \cong \R$ if and only if $K$ {\em compact}. 
It follows that for $K$ closed, connected and noncompact,
$\E(K)$ is nonempty. 

We may define ends of $K$ in geometric language as follows. A sequence
$(r_i)$ of properly embedded rays $r_i \subset P_X$ {\em limits to an end of $K$}
if for every open set $U$ containing $K$, there is a positive integer $N$ such
that for all $i,j \ge N$, the rays $r_i$ and $r_j$ are contained in $U$, and moreover for any
pair of proper sequences of points $p_k$ and $q_k$ in $r_i$ and $r_j$ there
is a proper (in $P_X$) sequence of arcs $\alpha_k \subset U$ joining $p_k$ to $q_k$.
Equivalently, for every {\em simply-connected} $U$ containing $K$, the rays
$r_i$ and $r_j$ should be properly homotopic in $U$ for all $i,j \ge N$.
Two such sequences limit to the same end, denoted $(r_i) \sim (r_i')$, if the
alternating sequence $r_1,r_1',r_2,r_2',\dots$ itself limits to an end of $K$.

\begin{defn}
For each $k \in T^\pm$, thinking of $k$ as a closed, connected
subset of $P_X$, define $\E_{k}$ to be the set of {\em ends} of $k$. Define 
$$\E^+ = \bigcup_{k \in T^+} \E_{k}$$
and define $\E^-$ similarly.
\end{defn}

\begin{lem}\label{R_tree_circular_ends}
There is a natural circular ordering on the set $\E^+$
which is preserved by the action of $\pi_1(M)$.
\end{lem}
\begin{proof}
Let $e_i \in \E_{k_i}$ be distinct elements of $\E^+$, for $i = 1,2,3$.
Of course, $k_i$ and $k_j$ might not be distinct if $e_i \ne e_j$, but
in this case, the ends $e_i,e_j$ of $k_i = k_j \subset P_X$ are different.

By removing a compact subset from $k_i$ in the case $k_i = k_j$, we can
replace the $k_i$ by disjoint connected sets $k_i'$ so that the $e_i$ are ends
of $k_i'$. Let $U(k_i')$ be disjoint open neighborhoods of the $k_i'$ in
$P_X$. Then we can find disjoint proper rays $r_i$ with $r_i \subset U(k_i')$
corresponding to the $e_i$.

We circularly
order the $r_i$ as follows. Let $D$ be a sufficiently large
closed disk in $P_X$ which intersects all the $r_i$. Let $p_i$ be the unique
points on $r_i \cap \partial D$ such that $r_i \backslash p_i$ consists of two
components, of which the unbounded one is disjoint from $D$. Then the
circular order on $\partial D$ defines a circular order on $p_i$ and therefore
on $r_i$.

To see that this is well-defined, suppose without loss of generality
that we replace $r_1$ with a different proper ray $r_1'$ corresponding to the same
end of $U(e_1)$. Then there are
a sequence of arcs $\alpha_j$ from $r_1$ to $r_1'$ which are contained
in $U(e_1)$, and which are proper in $P_X$. Such arcs are
disjoint from the other rays $r_2,r_3$, and therefore they certify that the
circular order of $r_1,r_2,r_3$ agrees with the circular order of $r_1',r_2,r_3$.
\end{proof}

\begin{lem}\label{stabilizer_cyclic}
Let $e \in \E^+$ be arbitrary. Then the stabilizer of $e$ in $\pi_1(M)$ is either
trivial or cyclic.
\end{lem}
\begin{proof}
If $\alpha \in \pi_1(M)$ fixes $e \in \E_{k}$, then $\rho_\hol(\alpha)$ must
fix $k$, and therefore $\alpha$ must fix $e^+(k) \in S^2_\infty$. But
since $M$ is a closed hyperbolic $3$--manifold, the stabilizer of any point
in $S^2_\infty$ is trivial or cyclic, as claimed.
\end{proof}

\begin{construct}\label{universal_circle}
We topologize the circularly ordered sets $\E^\pm$ with the order topology, and take their
completions $\overline{\E}^\pm$. 
This makes them into compact
circularly ordered sets, and therefore they are naturally order isomorphic to some
uncountable compact subset of $S^1$. 

By quotienting out complementary intervals to the image
of this subset, we get natural surjections to two circles
$$\overline{\E}^\pm \to (S^1_\u)^\pm$$
Since this construction is natural, the action of $\pi_1(M)$ on $\E^\pm$ extends
to an action on $\overline{\E}^\pm$ which factors through to an action on
$(S^1_\u)^\pm$.

The induced representations
$$\rho_\u^\pm\co \pi_1(M) \to \homeo((S^1_\u)^\pm)$$
are actually orientation preserving, since $\rho_\hol$ is orientation preserving.
\end{construct}

We now prove the first main result of the paper:

\begin{universal_circle_theorem}
Let $M$ be a closed oriented hyperbolic $3$--manifold with a quasigeodesic flow $X$.
Then there are faithful homomorphisms 
$$\rho_\u^\pm\co \pi_1(M) \to \homeo^+((S^1_\u)^\pm)$$
where $(S^1_\u)^\pm$ are topological circles,
called the {\em universal circles of $X$}.
\end{universal_circle_theorem}
\begin{proof}
In \fullref{universal_circle} we actually construct two
natural circles $(S^1_\u)^\pm$ and natural homomorphisms
$$\rho_\u^\pm\co  \pi_1(M) \to \homeo^+((S^1_\u)^\pm)$$
The maps $\E^\pm \to (S^1_\u)^\pm$ are $1$--$1$ outside of countably many
elements, so we can certainly find points in $(S^1_\u)^\pm$ in the
image of a unique element of $\E^\pm$. Call such ends {\em injective}. Moreover,
since the maps $e^\pm$ are nonconstant, the 
images $e^\pm(P_X)$ are uncountable.
It follows that there are points $p \in S^2_\infty$ in the image of $e^+$ which are
not fixed by any nontrivial $\alpha \in \pi_1(M)$, and for which every component $k$ of
$(e^+)^{-1}(p)$ has injective ends, and similarly for $e^-$. It follows that the 
natural actions of $\pi_1(M)$ on $(S^1_\u)^\pm$ are both faithful.
\end{proof}

\begin{cor}\label{Weeks_manifold}
The Weeks manifold does not admit a quasigeodesic flow.
\end{cor}
\begin{proof}
In \cite{Calegari_Dunfield}, It is shown that the fundamental group of the
Weeks manifold does not act faithfully on a circle. It follows from
\fullref{ucthm} that the Weeks manifold does not admit a quasigeodesic flow.
\end{proof}

\section{Properties of the universal circles}\label{sec5}

In this section, by studying the properties of $\sim^\pm$ in more detail, we
construct a pair of $\pi_1(M)$--invariant {\em laminations} (to be defined below)
$\Lambda^\pm$ for $(S^1_\u)^\pm$, and use this
structure to produce a single canonical circle $S^1_\u$ which maps
monotonically to each $(S^1_\u)^\pm$.

\subsection{Constructing laminations}

\begin{lem}\label{separating_class}
Some equivalence class $k$ of $\sim^+$ is separating in $P_X$.
\end{lem}
\begin{proof}
We suppose not and derive a contradiction. 

If each equivalence class $k$ is nonseparating, then $P_X\backslash k$ is
connected, and 
$\til{H}_0(P_X\backslash k)$ $= 0$. By \fullref{end_duality_theorem},
it follows that $\bar{H}^1_c(k) = 0$ and therefore $\bar{H}^0(k) \to \bar{H}^0_e(k)$
is surjective. Since $k$ is connected by hypothesis, and noncompact
by \fullref{preimages_noncompact}, this implies that $k$ has exactly one end.

Since each equivalence class $k$ of $\sim^+$ has a single end,
we can define a map $r\co  P_X \to (S^1_\u)^+$
by sending $k$ to this end. We show first that this map is continuous. Let
$k$ be an equivalence class, let $k^l,k^r$ be two other 
equivalence classes, and let $e,e^l,e^r$ denote the three unique ends of
these equivalence classes. We can join $k^l$ to $k^r$ by a compact arc
$\tau$ which avoids $k$. Let $I \subset (S^1_\u)^+$ be the open interval
complementary to $e^l,e^r$ containing $e$.

Suppose $p_i \to p \in k$ is a convergent sequence, and suppose $p_i$ is in the
equivalence class $k_i$. By the definition of the equivalence relation $\sim^+$, 
if $k_i \to K$ in the Hausdorff sense, then $k$ contains a connected component of $K$.
Moreover, the connected components of $K$ are contained in
equivalence classes of $\sim^+$.

If there is some other connected component $k'$, then either the unique ends of
$k'$ and $k$ are on the same side of 
$k^l \cup k^r \cup \tau$, or else $k_i$ intersects $\tau$
for sufficiently large $i$, and therefore $K \cap \tau$ is nonempty. In the second
case, since $\tau$ is arbitrary, it follows that $K$ intersects every arc
$\tau$ joining $k^l$ to $k^r$, and therefore some connected component of $K$
separates $k^l$ from $k^r$, contrary to our assumption that no equivalence class
of $\sim^+$ is separating.

It follows that the ends of $k'$ and $k$ are on the
same side of $k^l \cup k^r \cup \tau$, and therefore limit to points in $I$.
But $I$ was arbitrary, so $k$ and $k'$ limit to the same point in $(S^1_\u)^+$
and therefore the map $r$ is continuous.

We now show that there is a circle $\gamma \subset P_X$ such that the
map $r\co \gamma \to (S^1_\u)^+$ has degree $1$. Let $k_1,k_2,k_3$ be
three equivalence classes, so that their corresponding ends $e_1,e_2,e_3$
are distinct and positively ordered in $(S^1_\u)^+$.
Moreover, choose the $k_i$ such that the $e_i$ are the
unique ends corresponding to their images in $(S^1_\u)^+$. We can do this
since all but countably many ends are injective. Pick points $p_i \in k_i$
and join $p_i$ to $p_{i+1}$ by an arc $\alpha_i$ in the complement of
$k_{i+2}$ (here indices are taken mod $3$). 
Such arcs $\alpha_i$ exist, since each $k_i$ is nonseparating, by hypothesis.
Then define $\gamma = \alpha_1 \cup \alpha_2 \cup \alpha_3$. Since $\alpha_i$
is disjoint from $k_{i+2}$, the image $r(\alpha_i)$ avoids the point $r(k_{i+2})$.
It follows that $r\co \gamma \to (S^1_\u)^+$ has degree $1$.

But under any continuous map from $P_X$ to $S^1$, the image of any circle
in $P_X$ must map to $S^1$ by a degree $0$ map, since $P_X$ is contractible. 
This contradiction shows that some equivalence class is separating, as claimed.
\end{proof}

\begin{rmk}\label{compactification_remark}
We remark that there is a natural topology on $P_X \cup (S^1_\u)^+$
which gives it a compactification as a closed disk, and
similarly for $(S^1_\u)^-$. 

The topology is easy to define:
every {\em bounded} open set $U \subset P_X$ is still open.
If $p \in (S^1_\u)^+$ is arbitrary, we define a basis for the
topology near $p$ as follows. Let $I_i$ be a nested sequence of closed intervals
in $(S^1_\u)^+$ which converge to $p$ in the Hausdorff sense, and such
that the endpoints $(I_i)^l,(I_i)^r$ correspond to ends $e_i^l,e_i^r$ of
equivalence classes $k_i^l,k_i^r$ of $\sim^+$. Let $r_i^l,r_i^r$ be
proper rays contained in small open tubular neighborhoods of
$k_i^l,k_i^r$, which are eventually disjoint from every fixed equivalence class
$k$ of $\sim^+$ except for $k_i^l$ or $k_i^r$, and which
go out the ends corresponding to $e_i^l,e_i^r$. Let
$\tau_i$ be an arc joining the initial point of $r_i^l$ to the initial point
of $r_i^r$. Then the union $r_i^l \cup \tau_i \cup r_i^r$ separates $P_X$ into
two open sets $U,V$, one of which (say $U$) contains an end $e_{p'}$
of some $k_{p'}$ which corresponds to some $p' \in \inte(I_i)$.

Then we define $U \cup \inte(I_i)$ to be an open neighborhood of $p$
in $P_X \cup (S^1_\u)^+$. With this topology, the union is a Peano continuum;
that is, it is nondegenerate, perfectly separable and normal,
compact, connected and locally connected.
Now, a theorem of Zippin (see \cite[Theorem III.5.1]{Wilder}) characterizes the
$2$--disk in the following way:

\begin{thm}[Zippin]
Let $C$ be a Peano continuum containing a circle $J$ and satisfying the following
three conditions:
\begin{enumerate}
\item{$C$ contains an arc that spans $J$ (ie, intersects $J$ only at its endpoints)}
\item{Every arc of $C$ that spans $J$ separates $C$}
\item{No closed proper subset of an arc spanning $J$ separates $C$.}
\end{enumerate}
Then $C$ is homeomorphic to the closed $2$--disk with boundary $J$
\end{thm}

An arc $\overline{l}$ in $P_X \cup (S^1_\u)^+$ which spans $(S^1_\u)^+$ restricts to
a properly embedded line $l$ in $P_X$. The relevant separation
properties follow readily, and we can conclude that $P_X \cup (S^1_\u)^+$
is homeomorphic to a closed disk. This fact is not used elsewhere in this paper.
\end{rmk}

\begin{defn}
A {\em lamination} $\Lambda$ of $S^1$ is a closed subset of the
space of unordered distinct pairs of points in $S^1$ with the property
that for any two elements $\lbrace a,b\rbrace$ and $\lbrace c,d\rbrace$ 
of $\Lambda$, the pair of points $\lbrace a,b \rbrace$ does not link the pair of points
$\lbrace c,d \rbrace$ in $S^1$ (though it might share one or both points
in common). We also abbreviate this last condition by saying that $\Lambda$
is {\em unlinked} as a subset of the space of unordered distinct pairs of
points in $S^1$. We sometimes refer to elements of $\Lambda$ as {\em leaves}.
\end{defn}

\begin{rmk}
It is straightforward to define linking for pairs of 
disjoint $S^0$'s in any circularly ordered set.
\end{rmk}

The existence of a separating equivalence class lets us define naturally
a pair of laminations of $(S^1_\u)^\pm$.

\begin{construct}\label{circle_laminations}
Let $k$ be an equivalence class of $\sim^+$. If $k$ is separating, 
then $k$ has at least
two ends, as in the proof of \fullref{separating_class}. 
Let $\E_k$ denote the set of ends of $k$, and let
$\overline{\E_k}$ denote its closure in $\overline{\E}^+$. 

We let
$\Lambda_k$ be equal to the set of pairs of endpoints of closures of
complementary intervals of $\overline{\E_k}$.
Note that if $k \ne k'$ then $\E_k$ and $\E_{k'}$ do not link. For,
$k,k'$ are disjoint in $P_X$, so any two ends of $k'$ are contained on
the same side of $k$, and vice versa.
Since the stabilizer
of $k$ in $\pi_1(M)$ is at most cyclic, if we let $k$ be some separating
equivalence class of $\sim^+$, and $k' = \alpha(k)$ for some $\alpha$ not
in the stabilizer of $k$, we see that $\E_k$ is not dense in $\E^+$.

Conversely, if $U$ is a connected component of $P_X\backslash k$, 
then we can find $k' \subset U$
such that the image of $\E_{k'}$ is disjoint from $\E_k$ in $(S^1_\u)^+$.
Since we can do this for each component $U$, we see that distinct points of
$\overline{\E_k}$ cannot become identified under the quotient
map $\overline{\E}^+ \to (S^1_\u)^+$. It follows that
$\Lambda_k$ is nonempty as a lamination of $(S^1_\u)^+$.
$$\Lambda^+ = \overline{\bigcup_k \Lambda_{k}}
\leqno{\hbox{Now, define}}
$$
where $k$ ranges over equivalence classes of $\sim^+$.
\end{construct}

We summarize this construction in a lemma:

\begin{lem}\label{reps_fix_lamination}
The representations $\rho_\u^\pm$ preserve nonempty laminations
$\Lambda^\pm$ of $(S^1_\u)^\pm$.
\end{lem}
\begin{proof}
This follows from \fullref{separating_class} and
\fullref{circle_laminations}.
\end{proof}

\begin{lem}\label{compact_crossing}
Let $k^+,k^- \subset P_X$ be equivalence classes of $\sim^+,\sim^-$
respectively. Then the intersection $k^+ \cap k^-$ is {\em compact}.
\end{lem}
\begin{proof}
By definition, if $K = k^+ \cap k^-$ is the intersection, then every two
flowlines $l_1,l_2$ of $\til{X_\F}$ corresponding to points $l_1,l_2 \in K$
have the same endpoints in $S^2_\infty$. By \fullref{quasigeodesic_is_uniform}
and \fullref{straighten_estimate}, the flowlines $l_1$ and $l_2$ are 
distance $\le 2^9k^4$ apart in the Hausdorff metric. In particular, the
union $\til{X_\F}(K)$ of all flowlines corresponding to points in $K$
is itself a finite Hausdorff distance from a geodesic, and therefore
intersects a cross-section of the flow in a compact set.
\end{proof}

\fullref{compact_crossing} lets us directly compare the circles $(S^1_\u)^+$ and
$(S^1_\u)^-$.

\begin{lamination_theorem}
Let $M$ be a closed oriented hyperbolic $3$--manifold with a quasigeodesic flow $X$. Then
there is a canonical circle $S^1_\u$, a faithful homomorphism
$$\rho_\u\co \pi_1(M) \to \homeo^+(S^1_\u)$$
and natural monotone maps 
\begin{gather*}
\phi^\pm\co S^1_\u \to (S^1_\u)^\pm\\
\tag*{\hbox{so that}}
\phi^\pm_*(\rho_\u) = \rho_\u^\pm
\end{gather*}
Moreover, there are a pair of laminations $\Lambda^\pm_\u$ of $S^1_\u$
which are preserved by $\pi_1(M)$, and satisfy 
$$\phi^\pm(\Lambda^\pm_\u) = \Lambda^\pm$$
\end{lamination_theorem}
\begin{proof}
We define a circular order on $\E = \E^+ \cup \E^-$ as follows.
If $e_1,e_2,e_3$ are ends of $k_1,k_2,k_3$, each of which is an
equivalence class of {\em either}
$\sim^+$ or $\sim^-$, then by \fullref{compact_crossing}, outside
some compact ball $D$, the subsets $k_i$ are disjoint, and we
can take disjoint open neighborhoods around them and find proper rays
$r_i$ contained in these neighborhoods
which go out the ends corresponding to $e_i$.
Since the rays $r_i$ are disjoint, they admit a natural
circular order, which defines a circular order on $\E$, which restricts to
the circular order on each of $\E^+,\E^-$ obtained in \fullref{R_tree_circular_ends}.

As in \fullref{universal_circle}, 
we can take the order completion $\overline{\E}$ of $\E$, embed it in an
order-preserving way as a subset of $S^1$, and quotient out complementary
intervals to produce a circle $S^1_\u$.

The images of $\overline{\E}^\pm$ in $S^1_\u$ are closed subsets. If we further quotient
out complementary regions, we get quotient maps $\phi^\pm\co S^1_\u \to (S^1_\u)^\pm$.

It remains to define the laminations $\Lambda^\pm_\u$. For each leaf
$l = \lbrace a,b\rbrace$ of $\Lambda^+$ whose endpoints are in $\E^+$, 
we can include $\E^+$ into $\overline{\E}$ and map $\overline{\E}$ to $S^1_\u$,
and define $l_\u$ to be the leaf defined by the image of the endpoints.
Then define $\Lambda^+_\u$ to be the closure of the union of
the leaves $l_\u$ obtained this way.
Define $\Lambda^-_\u$ similarly.
\end{proof}

\subsection{Order trees}

To state the next corollary, we must recall the definition of an
{\em order tree}. 

\begin{defn}
An {\em order tree} is a set $T$ together with a collection $\mathcal{S}$ of
linearly ordered subsets called {\em segments}, each with distinct least
and greatest elements called the {\em initial} and {\em final} ends. If $\sigma$
is a segment, $-\sigma$ denotes the same subset with the reverse order, and
is called the {\em inverse} of $\sigma$. The following conditions should be
satisfied:
\begin{enumerate}
\item{If $\sigma \in \mathcal{S}$ then $-\sigma \in \mathcal{S}$.}
\item{Any closed subinterval of a segment is a segment (if it has more than one element).}
\item{Any two elements of $T$ can be joined by a finite sequence of
segments $\sigma_i$ with the final end of $\sigma_i$ equal to the initial end of
$\sigma_{i+1}$.}
\item{Given a cyclic word $\sigma_0\sigma_1 \cdots \sigma_{k-1}$ (subscripts
mod $k$) with the final end of $\sigma_i$ equal to the initial end of $\sigma_{i+1}$,
there is a subdivision of the $\sigma_i$ yielding a cyclic word 
$\rho_0\rho_1 \cdots \rho_{n-1}$ which becomes the trivial word when adjacent
inverse segments are cancelled.}
\item{If $\sigma_1$ and $\sigma_2$ are segments whose intersection is a single element
which is the final element of $\sigma_1$ and the initial element of $\sigma_2$
then $\sigma_1 \cup \sigma_2$ is a segment.}
\end{enumerate}
If all the segments are homeomorphic to subintervals of $\R$ with their order
topology, then $T$ is an {\em $\R$--order tree}.
\end{defn}

An order tree is topologized by the usual order topology on segments.
A lamination $\Lambda$ of a circle $S^1$ gives rise to a dual order tree
$\Lambda^*$ as follows: think of $S^1$ as the ideal boundary of $\H^2$, and
think of $\Lambda$ as the ideal points of leaves of 
a geodesic lamination $\Lambda_g$ of $\H^2$. The order tree $\Lambda^*$ is the
quotient of the space of leaves of $\Lambda_g$ by the relation which
identifies all boundary leaves of each complementary region. 
Each geodesic arc $\gamma \subset \H^2$
transverse to $\Lambda_g$ intersects $\Lambda_g$ in some closed subset, which inherits
an order structure from $\gamma$; this defines the set $\mathcal{S}$. Notice
that $\Lambda^*$ arising in this way is a {\em Hausdorff} order tree.

\begin{cor}
There are infinitely many closed hyperbolic $3$--manifolds without
quasigeodesic flows.
\end{cor}
\begin{proof}
In \cite{Fenley_noaction}, Fenley shows that there are infinitely
many hyperbolic surgeries $M_i$ on certain once-punctured torus 
bundles over the circle with the
property that every action of $\pi_1(M_i)$ on an order tree has a
global fixed point. Consider the manifold $M = M_i$. Suppose that
$M$ admits a quasigeodesic flow $X$. By \fullref{lamthm}, $\pi_1(M)$
acts faithfully on the order tree dual to the laminations
$\Lambda^+_\u$. But by Fenley, this action globally fixes
some point $p \in (\Lambda^+_\u)^*$.

It follows that under $\rho^+_\u$, the group
$\pi_1(M)$ fixes the union of the leaves $\lambda_i$ which are dual to $p$.
Such leaves correspond to the boundary leaves of the convex hull of the
set of ends of some equivalence class $k$ of $\sim^+$, and therefore
$\pi_1(M)$ must fix the corresponding point $e^+(k) \in S^2_\infty$.
But the stabilizer of any point in $S^2_\infty$ is cyclic or trivial; this
contradiction shows that $M$ does not admit a quasigeodesic flow, and
the theorem is proved.
\end{proof}

\section{Quasigeodesic flows and the Thurston norm}\label{sec6}

In this section we demonstrate that the Thurston norm on a hyperbolic
$3$--manifold $M$ can be characterized in terms of the homotopy classes of
quasigeodesic flows that $M$ supports. A basic reference for the Thurston
norm is \cite{Thurston_norm}. Another useful reference is \cite{Oertel}.

\subsection{Circle bundles, plane bundles, and Euler classes}

\begin{defn}
Let $X$ be a $C^1$ flow on an oriented $3$--manifold $M$. Let $\xi$ be a complementary
oriented $2$--plane field, so that $$TX \oplus \xi = TM$$ as oriented bundles.
Let $e(\xi) \in H^2(M;\Z)$ be the obstruction to trivializing $\xi$ as a bundle.
By abuse of notation, we denote the image of $e(\xi)$ in $H^2(M;\R)$ by $e(X)$.
\end{defn}

The class $e(\xi)$ is an example of a {\em characteristic class}, and is
called the {\em Euler class} of the bundle $\xi$.

Suppose $E$ is a bundle over a space $M$ with fiber homeomorphic to some fixed
manifold $F$. There is a classifying space $\bhomeo(F)$ for $F$ bundles, and such
bundles $E$ are classified by homotopy classes of maps $[M,\bhomeo(F)]$.
A cohomology class $e \in H^*(\bhomeo(F))$ pulls back under such a homotopy
class of map to a cohomology class on $M$, called a {\em characteristic class}.
The Euler class of a circle bundle is an example of such a characteristic class.
Both $\homeo^+(S^1)$ and $\homeo^+(\R^2)$ with the compact-open topology
are homotopic to $S^1$, and therefore 
$$H^*(\bhomeo^+(S^1)) = H^*(\bhomeo^+(\R^2)) = \Z[e]$$
where $e$ is free in dimension $2$. The pullback of this class is called
the {\em Euler class} of an oriented $S^1$ or $\R^2$ bundle over $M$.
Note that since $\bhomeo^+(S^1)$ and $\bhomeo^+(\R^2)$ are homotopic to
$\CP^\infty$ which is a $K(\Z,2)$, 
oriented plane or circle bundles are {\em classified} up
to isomorphism by the Euler class. Note further that if $M$ is a hyperbolic
$3$--manifold, $M$ is a $K(\pi,1)$, and therefore the Euler class may
be thought of as an element of group cohomology $H^2(\pi_1(M);\Z)$.
See \cite{Hus} or \cite{Milnor_Stasheff} for more details.

Associated to the representation $\rho_\u\co \pi_1(M) \to \homeo^+(S^1_\u)$
there is a (foliated) circle bundle $E_\u$ over $M$, defined by the
usual Borel construction
$$E_\u = \til{M} \times S^1_\u / 
(m,\theta) \sim (\alpha(m),\rho_\u(\alpha)(\theta))$$
where $\alpha$ ranges over $\pi_1(M)$, and $m \in \til{M},\theta \in S^1_\u$
are arbitrary. Similarly, the natural action of $\pi_1(M)$ on $P_X$ defines
a (foliated) $\R^2$ bundle $E_P$ over $M$, defined analogously.

For foliated $\R^2$ and $S^1$ bundles, we have the following algebraic
definitions of the Euler class, by thinking of elements of $H^2(G;\Z)$ as
central extensions:

\begin{construct}\label{circle_algebraic_class}
Let $G$ be a subgroup of $\homeo^+(S^1)$.
Let $\hat{G}$ denote the preimage of $G$ in the universal central
extension $\til{\homeo^+(S^1)}$ (the preimage of $\homeo^+(S^1)$ in $\homeo^+(\R)$ under
the universal covering map $\R \to S^1$). Then there is a central extension
$$0 \to \Z \to \hat{G} \to G \to 0$$
The class of this extension gives the 
Euler class $\rho^*[e] \in H^2(G;\Z)$.
\end{construct}

\begin{construct}\label{plane_algebraic_class}
Let $G$ be a subgroup of $\homeo^+(\R^2)$.
Let $G_\infty$ denote the {\em germ of $G$ at $\infty$}.
Let $A$ be a punctured disk neighborhood of $\infty$, and let $\til{A}$
be the universal cover of $A$. Then there is a central extension
$$0 \to \Z \to \hat{G}_\infty \to G_\infty \to 0$$
where $\hat{G}_\infty$ denotes the subgroup of periodic germs of homeomorphisms
of $\til{A}$ which cover elements of $G_\infty$. 
The class of this extension pulls back by the natural
restriction homomorphism $G \to G_\infty$ to give the 
Euler class $\rho^*[e] \in H^2(G;\Z)$.
\end{construct}

By considering the germs of equivalence classes of $\sim^\pm$ in $P_X$
at infinity, one sees that the Euler class of the foliated bundle $E_P$
is the obstruction to lifting the circular order on $\E^\pm$ to a
$\pi_1(M)$--invariant total order. It follows that the Euler classes
of $E_P$ and $E_\u$ are equal.

See \cite[Section 4.1]{Calegari_euler} for a justification of this construction, and for
more details.

Let $U\xi$ denote the circle bundle of unit vectors in $\xi$. Note that
$U\xi$ and $\xi$ have the same Euler class as $S^1$ and $\R^2$ bundles
respectively.

\begin{lem}\label{bundle_foliated}
The circle bundles $U\xi$ and $E_\u$ are isomorphic.
\end{lem}
\begin{proof}
At each point $x \in \til{M}$, the plane $\xi(x)$ may be thought of as the
tangent space at $x$ to the leaf space $P_X$ of $\til{X_\F}$. Note that
the hypothesis that $X$ is $C^1$ implies that $P_X$ admits a natural $C^1$
structure. By exponentiating, one sees that $\xi$ and $E_P$ are isomorphic as
plane bundles. 

Moreover, by the discussion above, the Euler classes of $E_P$ and $E_\u$
are equal. The claim follows.
\end{proof}

We will make use of a standard inequality known as the {\em Milnor--Wood inequality},
which is an inequality on the Euler class of a foliated circle bundle over
a surface. See \cite{Milnor} and \cite{Wood} for details.

\begin{thm}[Milnor--Wood]\label{Milnor--Wood}
Let $E$ be an oriented foliated circle bundle over a closed orientable surface
$S$. Then
$$|e(E)[S]| \le \max(0,-\chi(S)).$$
\end{thm}

\subsection{The Thurston norm}

With \fullref{Milnor--Wood} and \fullref{bundle_foliated} 
available to us, it is straightforward to prove the
main result of this section.

\begin{norm_theorem}
Let $M$ be a closed oriented hyperbolic $3$--manifold. Then the convex hull of
the set of Euler classes $e(X)$ as $X$ varies over the set of quasigeodesic
flows on $M$ is exactly the unit ball of the dual Thurston norm on $H^2(M)$.
\end{norm_theorem}
\begin{proof}
Let $X$ be a quasigeodesic flow on $M$, and let
$\xi$ be a complementary $2$--plane field. Then by \fullref{bundle_foliated},
$U\xi$ is isomorphic as a circle bundle to the foliated bundle $E_\u$.
Let $\Sigma$ be a surface in $M$ which is Thurston norm minimizing for
its homology class, so that $\|\Sigma\|_T = - \chi(\Sigma)$ where
$\|\cdot \|_T$ denotes the Thurston norm. The bundle $E_\u$ restricts to
a foliated circle bundle over $\Sigma$, so by the Milnor--Wood inequality
\fullref{Milnor--Wood}, we have $|e(X)[\Sigma]| \le \|\Sigma \|_T$;
that is, $e(X)$ is contained in the unit ball of the Thurston norm.

Conversely, Gabai and Mosher \cite{Mosher_flows}
constructed a pseudo-Anosov flow almost transverse
to any finite depth foliation. Fenley--Mosher \cite{Fenley_Mosher} showed
that these flows can all be taken to be quasigeodesic. Moreover, 
Gabai \cite{Gabai_finite_depth}
constructed a finite depth foliation containing any given norm-minimizing
surface as a compact leaf. This shows that for every $\Sigma$ as above,
there is {\em some} quasigeodesic flow $X$ such that
$$|e(X)[\Sigma]| = - \chi(\Sigma) = \|\Sigma\|_T.$$
Now, such $X$ are not necessarily $C^1$, but they can be $C^0$ approximated
by a $C^1$ flow $X'$. Such $X'$ is homotopic to $X$, and therefore has the
same Euler class. Moreover, since $k$ quasigeodesity is a local property
for any $k$, and since flowlines of $\til{X'_\F}$ are arbitrarily close
to flowlines of $\til{X_\F}$ on compact subsets, it follows that $X'$ is
also quasigeodesic.

In particular, every extremal point in the unit ball of the dual Thurston
norm can be realized, and therefore the conclusion follows.
\end{proof}

As remarked in the introduction, this generalizes one of the main theorems of
\cite{Mosher_efficient_intersection} and \cite{Mosher_norm_II}, for 
{\em pseudo-Anosov} quasigeodesic flows.

\bibliographystyle{gtart}
\bibliography{link}

\end{document}